
\documentstyle{amsppt}
\pagewidth{6in}
\vsize8.5in
\parindent=6mm
\parskip=3pt
\baselineskip=14pt
\tolerance=10000
\hbadness=500
\NoRunningHeads
\loadbold
\topmatter
\title
Singular Radon transforms and maximal  functions\\
 under
convexity assumptions\endtitle
\author
Andreas Seeger and Stephen Wainger
\endauthor
\thanks Research supported in part by  NSF grants.  \endthanks
\address
 Department of Mathematics, University of Wisconsin-Madison,
Madison, WI
53706, USA
\endaddress
\email seeger\@math.wisc.edu\endemail 
\address
 Department of Mathematics, University of Wisconsin-Madison,
Madison, WI
53706, USA
\endaddress
\email wainger\@math.wisc.edu \endemail
\abstract
We prove variable coefficient analogues of results in \cite{5} on Hilbert transforms and maximal functions  along convex curves  in the plane.
\endabstract
\endtopmatter

\define\inn#1#2{\langle#1,#2\rangle}

\define\lcontr{\rfloor}
\define\lco#1#2{{#1}\lcontr{#2}}
\define\lcoi#1#2{\imath({#1}){#2}}
\define\rco#1#2{{#1}\rcontr{#2}}

\define\lc{\lesssim}


\define\eps{\varepsilon}

\define\om{\omega}

\define\fM{{\frak M}}

\define\fa{{\frak a}}


\define\bbR{{\Bbb R}}

\define\bbZ{{\Bbb Z}}

\define\cC{{\Cal C}}

\define\cF{{\Cal F}}

\define\cH{{\Cal H}}

\define\cL{{\Cal L}}
\define\cM{{\Cal M}}

\define\cP{{\Cal P}}
\define\cQ{{\Cal Q}}
\define\cR{{\Cal R}}




\def\leaderfill{\leaders\hbox to 1em{\hss.\hss}\hfill}

\def\ga{{\gamma}}
\def\Ga{{\Gamma}}


\document
\head {\bf 1. Introduction}\endhead

The purpose of this paper is to prove $L^p$ boundedness results
on singular Radon transforms and maximal operators for variable curves in the plane.
We shall prove a diffeomorphism invariant extension of 
the result for translation invariant averages along along convex curves 
in \cite{5}.

To fix our notation  let $\Omega_0$, $\Omega_1$, $\Omega$  be  open sets in $\Bbb R^2$ with compact closure, so that $\Omega\subset\subset\Omega_1\subset\subset \Omega_0$. We assume that  for each 
$x\in \Omega_0$ we are given a curve
$$t\mapsto \Gamma(x,t), \qquad -c_0\le t\le c_0
\tag 1.1$$ so that
 $\Gamma(x,t)\in \Omega_0$ for all $x$ in   a
neighborhood of the closure of ${\Omega_1}$ and all $t\in [-c_0,c_0]$. 
Furthermore assume that
$\Ga$ satisfies
$$\Ga(x,0)=x,\tag 1.2$$
for all $x\in \Omega_0$. 
We denote by $\dot \Gamma(x,t)$ the $t$-derivative of $\Gamma$  and assume that $\dot \Gamma$ is an $L^\infty$ function, and 
that  $\Gamma$ and $\dot\Gamma$ depend smoothly on $x$.
We shall assume that for $|t|\le c_0$ the map $x\mapsto \Gamma(x,t)$ is a diffeomorphism on a neighborhood of $\Omega_1$ (for small $t$ this is of course implied by (1.2)). The inverse is denoted by $\Gamma^*$; thus
$x=\Gamma^*(y,t)$ iff $y=\Gamma(x,t)$.

The two  operators under consideration are  the maximal operator
$$\cM f(x)
=\sup_{0<h<\eps} \frac{1}{2h}\int^h_{-h}| f(\Gamma(x,t))| dt
\tag 1.3
$$
and the 
singular Radon transform
$$
\cR f(x)=\text{p.v.} \int \om(x,t)f(\Gamma (x,t)){dt\over t}
\tag 1.4
$$
where $\om$ is a $C^\infty_0$ function supported in  
$\Omega_0\times[-\eps,\eps]$. Here $\eps\le c_0$.
Since $\dot \Gamma$ is bounded it is not hard  to see that 
for $f\in C^1$
 the principal value integral (1.2) is well defined. 
Our task will be to show that under suitable assumptions the operators
$\cM$ and $\cR$ are $L^p$ bounded.
We  observe that it suffices to prove
$L^p$ estimates under the assumption that $\eps\ll c_0$ as the contribution for $t$ bounded  
away from $0$ is easy to handle.

As we are seeking to generalize the result in \cite{5} we 
 wish to make two assumptions on $\Gamma$,
namely a {\it convexity } hypothesis and a {\it doubling} hypothesis.
Since we consider
a variable 
 situation our assumptions ought to be invariant under changes of 
variables
(and the usual assumptions of convexity  fail to meet this requirement).

In order to introduce an  invariant convexity  assumption we 
follow \cite{23} and say 
that a function $h$ defined on an interval $J$ is {\it  quasi-monotonic} on $J$ if there is a constant
$\kappa\ge 0$  so that
$h'(t)=a(t)+E(t)$ for $t\in I$ where $a$ has constant sign in $I$ and $|E(t)|\le \kappa |h(t)|$
(typically $h$ is monotonic modulo a function in the ideal generated by $h$). A family of 
functions is {\it uniformly  quasi-monotonic} if in the latter inequality we can choose a 
universal $\kappa$.

The relevant quantities
are 
$$
\align
G(x,t)&=\det
\Big( \dot \Gamma(x,t)\quad\dot{\Gamma^*} (w,0)
\Big)
_{w=\Gamma(x,t)}
\tag 1.5
\\
G^*(y,t)&=\det\Big( \dot {\Gamma^*}(y,t)\quad\dot\Gamma (z,0)\Big)
_{z=\Gamma^*(y,t)}
\tag 1.6
\endalign
$$

We now make the following
\proclaim {Convexity Hypothesis (C.H.)}
For all $x\in \Omega_1$, $y\in  \Omega_1$
the functions  $G(x,\cdot)$ and $G^*(y,\cdot)$ are uniformly quasi-monotonic
on $[0,c_0]$ and on $[-c_0,0]$.
\endproclaim

We turn to  our doubling hypothesis.  We say that a non-negative  continuous 
function
$g$ on  $[0,c_0]$ is a {\it doubling function} if  $g(0)=0$, $g(t)>0$ for
$t>0$ and if there is $A\ge 1$ so that  $$g(t_2)\ge
2g(t_1) \quad\text{ if } t_2  \ge At_1.
\tag 1.7
$$
An immediate consequence is that
$$g(t_1)\lc (t_1/t_2)^\delta g(t_2), \qquad t_1\le A^{-1} t_2,\, t_2\le c_0, \tag 1.8 $$
for some $\delta>0$.

\proclaim{Doubling Hypothesis (D.H.)}
 There   is $C_0\ge 1$ and  a doubling function $g$ on $[0,c_0]$ so that
$$
C_0^{-1}g(A^{-1}|t|)\le |G(x,t)|\le C_0 g(A|t|)
\tag 1.9
$$
and
$$
C_0^{-1}g(A^{-1} |t|)\le |G^*(y,t)|\le C_0 g(A|t|).
\tag 1.10
$$
for all $x\in \Omega_1$,   $y\in \Omega_1$ and $|t|\le c_0$.
\endproclaim

In particular the inequality (1.8) 
holds for 
$G(x,\cdot)$ and $G^*(y,\cdot)$ if $t_1\le A^{-3} t_2 $, $t_2\le c_0$.

We can now formulate our main result.

\proclaim{Theorem A}  If the convexity hypothesis (C.H.) and the doubling
hypothesis (D.H.)  are satisfied then $\cM$ is bounded from  
$L^p$ to $L^p(\Omega)$, for 
$p>1$; moreover  $\cH$ is bounded from $L^p$ to $L^p(\Omega)$ 
 for $1<p<\infty$.
\endproclaim

Under very general finite type condition the $L^p$ boundedness of $\cM$ and $\cH$ has been 
proved by Christ, Nagel, Stein and Wainger  \cite{7} (see also Greenblatt \cite{11}). 
Thus we are mainly interested in the
flat case. The
 translation invariant model case of the 
theorem 
(where $\Gamma(x,t)=(x_1+t, x_2+u(t))$, with $u$ convex) 
was obtained  in Carlsson {\it et al.} \cite{5} ({\it cf.} also \cite{9}); the special case
$p=2$ goes back to \cite{19}, \cite{20} and in \cite{19} 
 it was also shown that our condition is necessary when
 $u$ is an even function. See also \cite{9} for a necessary condition in the general case. 
 In the `semi-translation invariant' case where
$\Gamma(x,t)=(x_1+t, x_2 +s(x_1,t))$ the $L^2$ result had been obtained by one of the authors 
in \cite{23}. 
$L^p$ theorems in somewhat different variable coefficient settings
are in \cite{3}, \cite{4}  and  in \cite{2}. More closely related to the setting here 
is the recent paper by Carbery and P\'erez \cite{1}
who proved $L^p$ bounds for the  semi-translation-invariant  case 
under   more restrictive third order assumptions.
Optimal results  on the Heisenberg group related to Theorem A
 were obtained 
by J. Kim \cite{15}, \cite{16}.

\subheading{Invariance properties and alternative formulations}
The main feature of hypotheses (C.H.) and (D.H.) is the invariance  under
diffeomorphisms. This is easy to check. Namely if $y=\Gamma(x,t)$, and
$y=\Phi(z)$, $x=\Phi(u)$, then $z=\widetilde \Gamma(u,t)$ with
$\widetilde \Gamma(u,t)=\Phi^{-1}(\Gamma(\Phi(u),t))$; moreover
$\widetilde \Gamma^*(u,t)=\Phi^{-1}(\Gamma^*(\Phi(u),t))$. Hence we get
$\dot{\widetilde \Gamma}(u,t)= D \Phi^{-1}_{\Gamma(\Phi(u),t)} \dot \Gamma(\Phi(u),t)$, and similarly 
$\dot{\widetilde \Gamma^*}(w,0)= D \Phi^{-1}_{\Gamma^*(\Phi(w),0)} 
\dot \Gamma^*(\Phi(w),0)$. The latter  we apply for $w=\widetilde \Ga(u,t)$ and notice that
$\Ga^*(\Phi(\widetilde\Ga(u,t)),0)=
\Ga^*(\Ga(\Phi(u),t),0)= \Ga(\Phi(u),t)$. Now let  $\widetilde G$ denote the determinant (1.5) corresponding to the 
curve $\widetilde \Ga$; then we obtain
$$\widetilde G(\Phi(u),t)=
\det(D\Phi^{-1}(\Gamma(\Phi(u),t)))\,
\det
\Big( \dot \Gamma(\Phi(u),t)\quad\dot{\Gamma^*} (\Ga(\Phi(u),t),0)\Big).
$$
A similar calculation applies to (1.6).  From this the invariance property easily follows,
with the possible change of the constants $A$, $C_0$  (see also the discussion in \cite{23}).

We also note the our assumptions do not depend on the particular parametrization. 
If $t=u(x,s)$ with $u_s\neq 0$ we have $\partial_s (\Gamma(x,u(x,s)))=u_s(x,s)\dot \Gamma(x,u(x,s))$ and the independence of the parametrization is easily verified.

Our  hypotheses 
 can also be described in terms of defining functions such as in
\cite{21}, \cite{24}. Namely let 
$\Sigma=\{(x,y): y=\Gamma(x,t), \text{ some } t\}$ then 
if we restrict to small values of $t$ 
the variety $\Sigma$ is a smooth hypersurface in $\Omega\times\Omega$ and 
$\Sigma=\{(x,y):\Psi(x,y)=0\}$ where $\Psi_x'\neq 0$ and $\Psi_y'\neq 0$.
Our  quasimonotonicity and doubling assumptions  
may be replaced by similar assumptions on the functions
$$
\align
&t\mapsto \text{det}(\Psi_y(x,y), \Psi_y(y,y))\Big|_{y=\Gamma (x,t)}
\tag 1.11
\\
&t \mapsto\text{det}(\Psi_x(x,y), \Psi_x (x,x))\Big|_{x=\Gamma^*(y,t)}. 
\tag 1.12
\endalign$$
If $N^*\Sigma\subset (T^*_L\Omega\setminus 0)\times (T^*_R\Omega\setminus 0)$ 
denotes the conormal bundle of $\Sigma$ then 
$$N^*\Sigma=\{ (x,\xi,y,\eta): \xi=\tau\Psi_x',\eta=\tau\Psi_y',  \tau\neq 0, \Psi(x,y)=0\}$$
and  assumptions on  (1.11-12) reflect  properties of the  
projections of $ N^*\Sigma$ to $ T_L^*\Omega$ and
$T_R^*\Omega$.

In order to see that the conditions involving (1.11-12) are equivalent to the conditions involving (1.5-6)
we first observe 
 that the conditions for (1.11-12) are invariant under changes of variables, moreover 
they do not depend on the particular choice of defining function.
By the above discussion we may  without loss of generality
 assume that 
$$\Gamma (x,t)=(x_1-t, \gamma(x_1, x_2,t)).
\tag 1.13
$$
  Then $\Gamma^*(y,t)=(y_1+t, \gamma^*(y_1, y_2, t))$ where
$\gamma (x_1, x_2, 0)=x_2$, $\gamma^*(y_1, y_2, 0)=y_2$,
${\partial\gamma\over \partial x_2}\not= 0$ and
${\partial\gamma^*\over\partial y_2}\not= 0$.  In fact
${\partial\gamma\over\partial x_2} (x_1, x_2,0)=1$ and
${\partial\gamma^*\over\partial y_2} (y_1, y_2, 0)=1$. 
The equivalence is now obtained by working with the defining functions 
$\Psi(x,y)= y_2-\gamma(x_1,x_2,x_1-y_1)$ or
$\widetilde \Psi(x,y)= x_2-\gamma^*(y_1,y_2,y_1-x_1)$.
These are both defining functions and they are related by
$$
y_2-\ga(x,x_1-y_1)=\fa(x,y)(x_2-\ga^*(y,y_1-x_1))
\tag 1.14
$$ where
$$ \fa(x,y)=
\int_0^1\frac{\partial\gamma}{\partial x_2}(x_1, (1-s)\gamma^*(y,y_1-x_1)+s x_2, x_1-y_1) ds.
$$
To see this expand $y_2-\gamma(x,x_1-y_1)$ about $x_2=\gamma^*(y,y_1-x_1)$ and use that
$$y_2=\ga(x_1,\ga^*(y,y_1-x_1),x_1-y_1).$$ Note that if $\eps$ is chosen small enough  we can 
assume that
$$|\fa(x,y)-1|\le 1/2 \qquad\text{ if } (x,y)\in\Omega_1\times \Omega_1, 
\ |x_1-y_1|\le \eps.
\tag 1.15 $$
For later reference we state that the boundedness of $\ga_{x_2}$ and $\nabla \ga_{x_2}$
(as assumed in Theorem B below)  imply that
$\fa$ has bounded derivatives.

\subheading{A change of variable}
The invariance under changes of variables allows us to 
to make a crucial choice of coordinates in order to reduce the situation  
(1.13) with the additional normalization
$\dot \gamma(x,0)=0$.
A related change  of coordinates was  suggested 
 years ago by C. Fefferman, in connection with the problem of differentiation 
along variable lines. A similar argument was also used in  \cite{25}.

We set  $\Phi (u_1, u_2)=(u_1,\rho(u_1,u_2))$ where the smooth function $\rho$ is
 to be determined and will satisfy
$\rho(0,u_2)=u_2$.   This also implies that for small $u_1$ the function
$u_2\mapsto \rho(u_1,u_2)$ is invertible, with inverse $\sigma$, so that
$\sigma(u_1, \rho(u_1,u_2))=u_2$.
Now suppose that were are already given $\Phi$ and  we would then have 
$$\Phi^{-1}\Gamma(\Phi(u),t)=(u_1-t, \sigma(u_1-t,\gamma (u_1, \rho(u_1,u_2), t)).
\tag1.16$$
Thus we 
 need to  take $\rho(\cdot, u_2)$ to satisfy the  equation
$$
-\sigma_{x_1}
(u_1,\gamma (u_1, \rho(u_1,u_2), 0))
+ \sigma_{x_2}(u_1,\gamma (u_1, \rho(u_1,u_2), 0)) \dot \gamma(u_1,\rho(u_1,u_2),0)=0.
\tag 1.17
$$
Now $\ga(u_1,\rho,0)=\rho$ and $\sigma_{x_1}(u_1,\rho(u_1,u_2)) +\sigma_{x_2}(u_1,\rho(u_1,u_2))\rho_{u_1}(u_1,u_2)=0$,
and thus (1.17) is implied by $\sigma_{x_2}\neq 0$ and
$$
\rho_{u_1}(u_1,u_2)
+ \dot \gamma(u_1,\rho(u_1,u_2),0)=0.
\tag 1.18
$$
Thus if we solve the ordinary differential equation (1.18), with parameter $u_2$, under the initial
 value condition $\rho(0,u_2)=u_2$ then we have $\rho_{u_2}\neq 0$ 
 and thus $\sigma_{x_2}(u_1,\rho(u))\neq 0$ for small $u_1$ and therefore
$\widetilde \gamma(u,t)=\sigma(u_1-t,\gamma (u_1, \rho(u_1,u_2), t))$ will satisfy 
$\dot{\widetilde \gamma}(u,0)=0$.

From now on we may and  shall work with 
families of curves defined by (1.13) which also satisfy
$$\dot \gamma(x,0)=0.\tag 1.19
$$
By implicit differentiation it also follows that
$$\dot{ \gamma^*}(y_1, y_2, 0)=0.
\tag 1.20$$

In this situation our convexity hypothesis simplifies to
$$
\align
\ddot{\gamma}(x, t)&=a(t,x)+ O(\dot \gamma(x,t)),\tag 1.21
\\
\ddot{\gamma^*}(y, t)&=a^*(t,y)+ O(\dot {\gamma^*}(y,t)),\tag 1.22
\endalign
$$
where $a(x,\cdot)$ and $a^*(y,\cdot)$ are  of constant sign for  $t>0$ 
and of constant sign for $t<0$.
Our doubling hypothesis becomes
$$\align
&C_0^{-1}g(A^{-1}t)\le |\dot\gamma(x, t)|
\le C_0g(A|t|),\tag 1.23
\\
&C_0^{-1}g(A^{-1}t)\le |\dot{\gamma^*}(y,t)|
\le C_0g(A|t|) \tag 1.24
\endalign
$$
for some doubling function $g$ and suitable constant $A\ge 1.$

We then have the following result:

\proclaim{Theorem B} Assume that $\gamma$ and $\gamma^*$ satisfy the
 hypotheses (1.19-24).  Also suppose $|\partial_{x_2}\gamma (x,t)|\ge c_1>0$.
Then $\cM$ is bounded from  
$L^p$ to $L^p(\Omega)$, for 
$p>1$, and   $\cH$ is bounded from $L^p$ to $L^p(\Omega)$. The operator norms
depend only on the cutoff function $\om$,
the doubling function $g$, the constants $A$ in (1.23-24)
and the $L^\infty$ norms of $\gamma_{x_2}$, $\nabla\gamma_{x_2}$,
 $\gamma^*_{y_2}$, $\nabla\gamma^*_{y_2}$.
\endproclaim

With the change of variables discussed above, Theorem B implies Theorem A.

\remark{ Remark}  Note that the operator norms do not 
explicitly depend on the $L^1$ norm of $\ddot \gamma$. Thus by limiting arguments
 Theorem $B$ covers examples such as
 $\gamma(x,t)=x_2+u(t)$ where $u$ is even or odd, continuous, linear on $(2^{-j},2^{-j+1})$  with 
$u(2^{-j})=2^{-m j} $, as well as variable perturbations.
\endremark

\medskip

The organization of the paper is as follows.  In Section 2 we introduce
some notation and make a preliminary Littlewood Paley decomposition of our
operators; moreover we prove the $L^p$ estimates for the `Calder\'on-Zygmund part' of the operator.  In section 3 we give an outline of the proof of Theorem B, and handle the technical details of the main error estimate in \S4.

\head{\bf 2. Preliminary decompositions and Calder\'on-Zygmund estimates}
\endhead

Let $\phi\in C^\infty_0(\Bbb R)$ be supported in $(1/2,2)\cup(-2,-1/2)$
and define, $\phi_j(s)= 2^j\phi(2^js)$ for $j>0$.

Also let $\chi\in C^\infty_0(\Omega\times\Omega)$,
$\widetilde \chi$, $\widetilde \phi$ nonnegative so that $0\le |\chi|\le \widetilde \chi$,
$0\le |\phi|\le \widetilde \phi$ and 
$\widetilde \phi_j= 2^j\widetilde \phi(2^j\cdot)$. Define
$$
\aligned
\cR_j f(x) &=\int \phi_j(x_1-y_1) \chi (x,y_1, x_2+\gamma(x,x_1-y_1))
f(y_1, x_2+\gamma(x,x_1-y_1)) dy_1
\\
\cM_j f(x) &=\int \widetilde \phi_j(x_1-y_1) \widetilde \chi (x,y_1, x_2+\gamma(x,x_1-y_1))
f(y_1, x_2+\gamma(x,x_1-y_1)) dy_1;
\endaligned
\tag 2.1
$$
here  we want to explicitly include the case that the functions $\widetilde \chi$ and $\chi$ and the functions
 $\widetilde \phi$ and $\phi$  coincide (and are nonnegative).

The $L^p$ inequality for the maximal function in (1.3)  is a simple 
consequence  of the $L^p$ boundedness of the maximal operator $\cM$ defined with slight abuse of notation by
$$\cM f(x)=\sup_{j\in J}| \cM_j f(x)|;$$
here
$J$ is a finite set of  integers $j>C$
 and the bound is not supposed to depend on the cardinality of $J$. 
Working with suitable positive cutoff functions  we obtain uniform bounds for $\cM$ from uniform bounds for the maximal function 
$$\sup_{j\in J} |\cR_j f|.$$ 
Notice that since every individual operator 
$\cR_j$ is bounded on $L^1$ and $L^\infty$ we 
need to take the supremum over large $j$ only.
Similarly the boundedness of the Hilbert transform follows from 
 the $L^p$ boundedness of the
operator $\sum_{j>C} \cR_j f$ under the additional assumption
that the cutoff function satisfies $\int \phi(s)=0$; indeed we can choose 
 $\phi$ such that $\sum_{j=-\infty}^\infty \phi_j(s)=1/s$.

Denote by $\delta_0$ the Dirac measure on the real line, at the origin.
Following \cite{21} we express $\delta_0(y_2-\gamma(x,t))$  as an 
oscillatory integral distribution
using the 
Fourier inversion formula,
$$\delta_0(y_2-\ga(x,t))=(2\pi)^{-1}
\int e^{i\tau (y_2-\ga(x,t))} d\tau,$$
 and then   decompose the singular integral operator  as in \cite{23} 
into two parts,
a low frequency part where the cancellation of $\phi$ is crucially  
used, and a high frequency part where this cancellation 
does not play a role.
See also \cite{18}, \cite{12}, \cite{22} for earlier 
variants of this approach.
The analogous decomposition is made for the maximal operator 
where of course no cancellation of $\phi$ is needed.

In order to proceed with this decomposition  we set $B=2^{20} A$ where $A$ is the constant in (1.7-10)
and  define integers $a_j$, $b_j$ so that
$$
\aligned
&2^{-a_j-1}<2^{-j}g(2^{-5-j} B)\le 2^{-a_j}
\\
&2^{-b_j-1}<2^{-j}g(2^{5-j} B^{-1})\le 2^{-b_j}
\endaligned
\tag 2.2
$$

For later reference we note that  for   $2^{j-k}\le (4A)^{-1}$ we have
$2^{a_j-a_k}\lc 2^{j-k}$; this does not use the full strength of the doubling assumption as it follows from (1.8) with  $\delta=0$. By the doubling assumption  the former estimate can be improved to
$2^{a_j-a_k}\lc 2^{(j-k)(1+\delta)}$, for some $\delta>0$.

Let $\beta_0$ be an even function in $C^\infty_0(\bbR)$ so that $\beta_0(s)=1$ if $|s|\le 1/2$ 
and  $\beta_0(s)=0$ if $|s|\ge 3/4$ and let 
for $\beta_k(s)=\beta_0(2^{-k}s)-\beta_0(2^{-k+1}s)$, for $k\ge 1$, 
so that $\beta_k(s)$   can be  nonzero for $k\ge 1$  only when $2^{k-1}< s < 2^{k+1}$. Clearly we  have $\sum_{k=0}^\infty \beta_k(s)\equiv 1$.

For $k>1$ we define operators 
$\cR_j^k$ with distribution kernel
$$
\align
R_j^k(x,y)&=\chi(x,y) \phi_j(x_1-y_1)\int e^{i\tau (y_2-\ga(x,x_1-y_1))}\beta_k(\tau) d\tau
\endalign
$$
as well as operators $\cH_j$ with distribution kernel
$$H_j(x,y)=
 \chi(x,y) \phi_j(x_1-y_1)
\int e^{i\tau (y_2-\ga(x,x_1-y_1))}\beta_0(2^{-a_j}\tau) d\tau;
$$
then our basic decomposition is given by
$$2\pi \cR_j= \cH_j+ \sum_{k>a_j}\cR_j^k
\tag 2.3
$$
In the remainder of this section we shall first prepare further 
the term $\cR_j^k$ which for $k>a_j$ will later be treated as a piece of singular Fourier integral operator  and then deal with the
contribution  $\sum_j \cH_j$ (or the associated maximal function) 
which corresponds to a kind of
Calder\'on-Zygmund operator.

In (2.3) the decomposition in $k$ corresponds essentially to a Littlewood-Paley decomposition in the variable dual to $x_2$. To make this precise 
we introduce a Littlewood-Paley operator  $\cL^k$ defined by
$$
\widehat {\cL^k f}(\xi) =
\bigl[\beta_0(2^{-k-10}\xi_2)-\beta_0(2^{-k+10}\xi_2)\bigr] \widehat f(\xi)
\tag 2.4
$$
so that the multiplier is supported where $2^{k-10}\le |\xi_2|\le 2^{k+11}$ and equals $1$
on $2^{k-9}\le |\xi_2|\le 2^{k+9}.$

\proclaim{Lemma 2.1} For $1\le p\le \infty$
$$\|\cR_j^k- \cL^k\cR_j^k\cL^k\|_{L^p\to L^p}\lc 2^{-k}
$$
\endproclaim

Lemma  2.1 tells us that for $k\ge a_j$ we may replace the operators
$\cR_j^k$ by $\cL^k \cR_j^k \cL^k$ since
$$
\sum_{j}\sum_{k\ge a_j}
\|\cR_j^k- \cL^k\cR_j^k\cL^k\|_{L^p\to L^p}\lc
\sum_{j}\sum_{k\ge a_j} 2^{-k}
\lc \sum_j 2^{-a_j}\lc 1.
\tag 2.5
$$

\demo{\bf Proof of Lemma 2.1} We write 
$\cR_j^k- \cL^k\cR_j^k\cL^k= \cR^k_j(I-\cL^k) +(I-\cL^k)\cR^k_j \cL^k$. Thus it suffices to 
show that 
$ \cR^k_j(I-\cL^k)$ and  $(I-\cL^k)\cR^k_j$ satisfy the asserted bounds.  Let $\cP^{l,2}$ be the convolution operator with Fourier multiplier 
$\beta_0(2^{-l}\xi_2)$ and let
$\cQ^{l,2}$ be the convolution operator with 
 multiplier $\beta_1(2^{-l}\xi_2)$. Then  by the support properties 
the symbols we have
$$I-\cL^k= (I-\cL^k)(\cP^{k-5,2}+\sum_{l\ge k+ 5} \cQ^{l,2})$$
and consequently it suffices to show that
$$
\align
&\|\cR^k_j \cQ^{l,2}\|_{L^p\to L^p}
+ \|\cQ^{l,2}\cR^k_j \|_{L^p\to L^p}
\lc 2^{-l} \text{ if } l\ge k+5
\tag 2.6.1
\\
&\|\cR^k_j \cP^{k-5,2}\|_{L^p\to L^p}
+ \|\cP^{k-5,2}\cR^k_j \|_{L^p\to L^p}
\lc 2^{-k}
\tag 2.6.2
\endalign
$$
These estimates follow by standard  integration by parts arguments (see \cite{13}). For the sake of
 completeness we include the argument.
We first estimate the kernel 
of the operator 
$\cR^k_j \cQ^{l,2}$ which is given by
$$
K_{kl}(x,z)= \phi_j(x_1-z_1)
\iiint \chi(x,z_1,y_2)  e^{i(\tau(y_2-\gamma(x,x_1-z_1))+\xi_2(y_2-z_2))}
\beta_k(\tau)\beta_l(\xi_2) dy_2 d\tau d\xi_2 .
$$
Note that on the support of the symbol we have $|\xi_2+\tau|\approx 2^l$ . We integrate by parts 
once with respect to $y_2$ and then we integrate by parts with respect to $\tau$ and $\xi_2$. This yields the bound
$$
|K_{kl}(x,z)|\lc|\phi_j(x_1-z_1)| 2^{-l} 
\int\frac{2^k}{(1+2^k|y_2-\gamma(x,x_1-z_1)|)^N}
\frac{2^l}{(1+2^l|y_2-z_2|)^N} dy_2
$$
and integration with respect to $z$ yields that $\int|K_{kl}(x,z)|dz\lc 2^{-l}$ uniformly in $x$.
If we take into account (1.15)
 then we also get 
that
$\int|K_{kl}(x,z)|dx\lc 2^{-l}$ uniformly in $z$ and the asserted bound (2.6.1) for 
$\cR^k_j \cQ^{l,2} $ follows. The proof of (2.6.2) for
$\cR^k_j \cP^{k-5,2}$ is the same.

Next we examine the kernel of
the operator 
$\cQ^{l,2} \cR^k_j $
which is given by
$$
\widetilde K_{kl}(x,z)=
\phi_j(x_1-z_1) 
\iiint \chi(x_1,w_2,z)  
e^{i(\xi_2(x_2-w_2)+\tau(z_2-\ga(x_1,w_2,x_1-z_1)))}
\beta_k(\tau)\beta_l(\xi_2) dw_2 d\tau d\xi_2 .
$$
The difference is now the nonlinear dependence on the phase in $w_2$. To remove this potential difficulty we may again  invoke (1.15)
 and change variables
in the oscillatory integral to $\sigma=\tau \fa(x_1,w_2, z)$. Thus we get
$$\multline
\widetilde K_{kl}(x,z)=\\
\phi_j(x_1-z_1) 
\iiint \frac{\chi(x_1,w_2,z)} {\fa(x_1,w_2,z)}   
e^{i(\xi_2(x_2-w_2)+\sigma(w_2-\ga^*(z,x_1-z_1)))}
\beta_k(\frac{\sigma}{\fa(x_1,w_2,z)})\beta_l(\xi_2) dw_2 d\tau d\xi_2 .
\endmultline
$$
With this representation the estimation of 
$\widetilde K_{kl}$ is exactly the same as for $K_{kl}$.
Recall  that $|\fa-1|\le 1/2$. In the integration by parts 
with respect to $w_2$ we shall also need the boundedness of 
$\partial_{x_2}\fa$ which is guaranteed by our assumption,
{\it cf.} the remark following (1.15).
As above we see that the $L^p\to L^p$ bound for 
$\cR^k_j \cP^{k-5,2}$ is $O(2^{-k})$. The proof of the bound
$\|\cP^{k-5,2}\cR^k_j \|_{L^p\to L^p}=O(2^{-k})$ is the same.\qed
\enddemo

Concerning the operators $\cH_j$ we make the following simple observation (which is valid 
without any cancellation property). 

\proclaim{Lemma 2.2} The kernels $H_j$ satisfy 

$$\align
& |\partial_{y_1}^{n_1}
\partial_{y_2}^{n_2}H_j(x,y)\big| 
 \le A_0
2^{jn_1+a_jn_2}  2^{j+a_j}(1+2^j|x_1-y_1|)^{-2}
(1+2^{a_j}|x_2-y_2|)^{-2},\tag 2.7.1
\\
& |\partial_{x_1}^{n_1}
\partial_{x_2}^{n_2}H_j(x,y)\big| 
 \le A_0 
2^{jn_1+a_jn_2}  2^{j+a_j}(1+2^j|x_1-y_1|)^{-2}
(1+2^{a_j}|x_2-y_2|)^{-2},\tag 2.7.2
 \endalign
$$
 for $(n_1,n_2)=(0,0)$, $(0,1)$ or $(1,0)$.
\endproclaim

\demo{\bf Proof} 
By integration by parts we have 
$$|H_j(x,y)| \lc 2^j\chi_{[-2^{1-j},2^{-j+1}]}(x_1-y_1)
2^{a_j}(1+2^{a_j}|y_2-\gamma(x,x_1-y_1)|)^{-N} $$
Now  if $|x_1-y_1|\le 2^{-j+1} $ then 
 $$\align
|\gamma(x, x_1-y_1)-x_2|&=|x_1-y_1|\Big|\int_0^1 \dot \gamma(x, s(x_1-y_1)) ds\Big|
\\&\lc |x_1-y_1| g(A2^{-j+1}) \lc 2^{-a_j}
\endalign
$$
and thus 
$1+2^{a_j}|y_2-\gamma(x,x_1-y_1)|\approx
1+2^{a_j}|y_2-x_2|$ if  $|x_1-y_1|\le 2^{-j+1}$. This yields the asserted estimate for 
$n_1=n_2=0$.
The estimates for the derivatives are analogous.\qed
\enddemo

Let $M_{str}$ be the strong maximal operator 
(involving averages over rectangles parallel to the coordinate axes).
Then the following
 estimate is an immediate consequence of Lemma 2.2.

\proclaim{Corollary  2.3} For all $x\in \Omega$,
$$\sup_j |\cH_j f(x)|\lc M_{str} f(x).\tag 2.8 $$ 
\endproclaim

In the case where $\int \phi(s) ds=0$ we get a bound for the sum $\sum_j \cH_j$. As in  \cite{2} 
the $L^p$ boundedness
is proved by invoking the  Hardy space $H^1_{prd}:=H^1(\bbR\times \bbR)$  defined using the
two-parameter dilations (see \cite{6}, \cite{14}, \cite{10}). 
Recall from \cite{6} that operators which are bounded on $H^1_{prd}$ and bounded on 
$L^2$ are also bounded on $L^p$ for $1<p<2$.

\proclaim{Proposition 2.4} Suppose that the cancellation condition 
 $\int \phi(s) ds =0$ holds. Then  the
operators $\sum_j \cH_j$ and
 $\sum_j \cH_j^*$ are both bounded on $L^2$ and on $H^1_{prd}$, 
and consequently on
 $L^p$, for $1<p<\infty$.
\endproclaim

\demo{\bf Proof} We first show the  $L^2$ boundedness. By the almost orthogonality 
lemma of Cotlar and Stein  it suffices to verify
$$\|\cH_j \cH_k^*\|_{L^2\to L^2}
+\|\cH_j^* \cH_k\|_{L^2\to L^2}
\lc 2^{-|j-k| /2}.
\tag 2.9$$
By taking adjoints it suffices to show (2.9) for $k\ge j$ and since the operator norms of $\cH_j$ are uniformly bounded it suffices to consider the case where $2^{k-j}\ge 2^{10}A$.

We first examine
$\cH_j \cH_k^*$; its kernel is given by 
$$
\int H_j(x,z)\overline {H_k(y,z)} dz 
=\sum_{n\ge 0} H_{jk}^n(x,y)
$$
where
$$\align  H_{jk}^n(x,y)
=\iiint
\phi_j(x_1-z_1) \phi_k(y_1-z_1) 
e^{i( \tau(z_2-\ga(x,x_1-z_1))-\sigma (z_2-\ga(y,y_1-z_1)))}&
\\
\\ \times \
\beta_0(2^{-a_j}\tau)
\beta_0(2^{-a_k}\sigma) \beta_n(2^{-a_j}\sigma)
dz d\tau d\sigma.&
\tag 2.10
\endalign
$$
Here we used that 
$\sum_{n=0}^\infty \beta_n(2^{-a_j}\sigma)\equiv 1$. Observe that 
in view of the support properties of the symbol
we have the restriction $a_j+n\le a_k+1$.

Now let
$h(z_1)=
\phi_j(x_1-z_1)
e^{i( \tau(z_2-\ga(x,x_1-z_1))-\sigma (z_2-\ga(y,y_1-z_1)))}$. We use the cancellation of 
$\phi_k$ to replace 
$h(z_1)$ in (2.10) by $h(z_1)-h(y_1)=(z_1-y_1)\int h'(y_1+s(z_1-y_1)) 
ds$; this will be relevant for small $n$.
We write $\zeta_{j,k,n}(\sigma,\tau)=\beta_0(2^{-a_j}\tau)
\beta_0(2^{-a_k}\sigma) \beta_n(2^{-a_j}\sigma)$ and obtain
$$ H_{jk}^n(x,y) = \int_0^1 [I_s(x,y)+II_s(x,y)+ III_s(x,y)] ds
$$
where
$$
\align
I_s(x,y) &=\iiint
e^{i( \tau(z_2-\ga(x,x_1-z_1))-\sigma (z_2-\ga(y,y_1-z_1)))}
\zeta_{j,k,n}(\sigma,\tau) \chi_{j,k,1}(x,y,s)
dz d\tau d\sigma,
\\
II_s(x,y) &=\iiint
e^{i( \tau(z_2-\ga(x,x_1-y_1+s(y_1-z_1)))-\sigma (z_2-\ga(y,y_1-z_1)))}
 \zeta_{j,k,n}(\sigma,\tau)
\tau \chi_{j,k,2}(x,y,s)
dz d\tau d\sigma,
\\
III_s(x,y) &=\iiint
e^{i( \tau(z_2-\ga(x,x_1-z_1))-\sigma (z_2-\ga(y,s(y_1-z_1))))}
\zeta_{j,k,n}(\sigma,\tau)\sigma  \chi_{j,k,3}(x,y,s)
dz d\tau d\sigma.
\endalign
$$
with
$$
\align
\chi_{j,k,1}(x,y,s)&=
(y_1-z_1) \phi_j'(x_1-y_1+s(y_1-z_1)) \phi_k(y_1-z_1) 
\\
\chi_{j,k,2}(x,y,s)&=
\phi_j(x_1-z_1) \phi_k(y_1-z_1)  (y_1-z_1) \dot\gamma(x,x_1-y_1+s(y_1-z_1))
\\
\chi_{j,k,3}(x,y,s)&=
\phi_j(x_1-z_1) \phi_k(y_1-z_1)  
(y_1-z_1) \dot\gamma(y,s(y_1-z_1)).
\endalign
$$
Since $|y_1-z_1|\le 2^{-k+1}$  we then have
$$
\align
&|\chi_{j,k,1}(x,y,s)|\le 2^{2j-k}
\\
&
2^{a_j} |\chi_{j,k,2}(x,y,s)|
\lc 2^{a_j} g(2^{1-j}A)2^{-k} \lc  2^{j-k}
\\
&2^{a_j+n}|\chi_{j,k,3}(x,y,s)|
\lc 2^{a_j+n} g(2^{1-k}A)2^{-k} \lc 2^{a_j-a_k+n}.
\endalign
$$

These estimates are used after  additional integration by parts in $\tau$ and
 $\sigma$. For the term $I_s$ we obtain (taking into account the symbol properties of $\zeta_{j,k,n}$)

$$\multline |I(x,y)|\lc 2^{2j-k} 2^{k}
\\ \times
\iint\limits\Sb|x_1-z_1|\le 2^{-j+1}\\|y_1-z_1|\le 2^{-k+1}\endSb
\frac{2^{a_j}}{(1+2^{a_j}|z_2-\gamma(x,x_1-z_1)|)^N}
\frac{2^{a_j+n}}{(1+2^{a_j+n}|z_2-\gamma(y,y_1-z_1)|)^N}
dz_1dz_2.
\endmultline
$$
Observe that the 
integral $\int_{x_2} 2^{a_j}(1+2^{a_j}|z_2-\gamma(x,x_1-z_1)|)^{-N} dx_2$ 
is $O(1)$ in view of (1.15).
Thus in evaluating
$\int |I(x,y)| dx$, for fixed $y$, we perform an $x_2$ integration first and see that
$$\align
\int|I(x,y)| dx&\lc
2^{2j}
\iint\limits\Sb|x_1-z_1|\le 2^{-j+1}\\|y_1-z_1|\le 2^{-k+1}\endSb
\int 
\frac{2^{a_j+n}}{(1+2^{a_j+n}|z_2-\gamma(y,y_1-z_1)|)^N} dz_2 \, dz_1 dx_1\lc 2^{j-k}
\endalign
$$

We argue similarly for the terms $II_s$ and $III_s$.
By integration by parts we get the pointwise estimate
$$\multline |II_s(x,y)|\lc 2^{j+k} 2^{a_j-a_k}
\\ \times
\iint\limits\Sb|x_1-z_1|\le 2^{-j+1}\\|y_1-z_1|\le 2^{-k+1}\endSb
\frac{2^{a_j}}{(1+2^{a_j}|z_2-\gamma(x,x_1-z_1+s(z_1-y_1))|)^N}
\frac{2^{a_j+n}}{(1+2^{a_j+n}|z_2-\gamma(y,y_1-z_1)|)^N}
dz_1dz_2
\endmultline
$$
and 
$$\multline |III_s(x,y)|\lc 2^{j+k} 2^{a_j-a_k+n}
\\ \times
\iint\limits\Sb|x_1-z_1|\le 2^{-j+1}\\|y_1-z_1|\le 2^{-k+1}\endSb
\frac{2^{a_j}}{(1+2^{a_j}|z_2-\gamma(x,x_1-z_1)|)^N}
\frac{2^{a_j+n}}{(1+2^{a_j+n}|z_2-\gamma(y,s(y_1-z_1))|)^N}
dz_1dz_2.
\endmultline
$$
Since $2^{a_j-a_k}\lc 2^{j-k}$  we obtain the same bound 
$O(2^{j-k})$ for $\int |II_s(x,y)|dx$ as above, similarly for 
$\int |III_s(x,y)|dx$ we obtain the bound $O(2^{a_j-a_k+n})$ 
which is $O(2^{j-k+n})$. Thus 
$$
\int |H^n_{jk}(x,y)| dx\le \int_0^1\int|I_s|+|II_s|+|III_s| dx ds \lc 2^{j-k+n}.
\tag 2.11$$
The same bound is obtained for 
$\int |H^n_{jk}(x,y)| dy$, uniformly in $x$.

For large $n$ the estimate (2.11) is not sufficient 
but we can now use an integration by parts 
in $z_2$, in order to gain a factor   $2^{-a_j-n}$; this is followed as above by integration by 
parts with respect to $\tau$ and $\sigma$. The result is 
that for $n\ge 10$
$$
\int|H_{jk}^n(x,y)| dx\lc 2^{-a_j-n}
$$
uniformly in $y$ and again the same bound holds also for
$\int |H_{jk}^n(x,y)| dy$, uniformly in $x$.

We sum in $n$ and obtain the bound
$$\|\cH_j \cH_k^*\|_{L^2\to L^2} \le C \Big(
2^{j-k}+2^{a_j-a_k}
+\sum_{n\ge 10} \min\{ 2^{-a_j-n}, 2^{j-k+n}\}\Big).
$$
Now  if $n\ge 10$ we  use the bound  $2^{-a_j-n}$ for $n>  (k-j)/2$ and the bound
$2^{j-k+n}$ for $n\le (k-j)/2$. We sum in $n$ and
obtain  the asserted bound (2.9) for  the term
$\|\cH_j \cH_k^*\|$.

The estimation of $\|\cH_j^* \cH_k\|$ is largely analogous. However we first 
use (1.15)
 to represent the kernel of $\cH_j$ as
$$H_j(v,w)= \phi_j(v_1-w_1) \chi(v,w) \int 
e^{i\tau(v_2-\gamma^*(w,w_1-v_1))}
\beta_0(2^{-a_j}\frac{\tau}{\fa(v,w)})\frac{d\tau}{\fa(v,w)}.
$$ 
Thus the kernel of $\cH_j
^* \cH_k$ is given by
$$\align
\int &\overline{H_j(z,x)}H_k(z,y) dz=
\iiint \overline{\phi_j(z_1-x_1)\chi(z,x)}
\phi_k(z_1-y_1)\chi(z,y)\times
\\
&e^{-i(\tau(z_2-\ga^*(x,x_1-z_1))
-\sigma(z_2-\ga^*(y,y_1-z_1)))}
\overline{\beta_0(2^{-a_j}\frac{\tau}{\fa(z,x)})}
\beta_0(2^{-a_k}\frac{\sigma}{\fa(z,y)})
\frac{d\tau}{\fa(z,x)}
\frac{d\sigma}{\fa(z,y)} dz
\endalign
$$
and by  using this expression the above proof for
$\cH_j \cH_k^*$ can be repeated here. Again the only difference  is that we have to take 
into account  the limited differentiability of the symbol, 
but our assumptions on $\gamma_{x_2}$ and its gradient 
still allow us to once  integrate by parts with respect to $z_2$.\qed
\enddemo

In order to complete the proof of the $H^1_{prd}\to L^1$ boundedness we use the following
Lemma which is  proved by standard arguments.

\proclaim{Lemma 2.5} Let $\{H_j\}_{j\in I}$ be a finite family of Schwartz kernels
and let $\cH_j$ be the associated operators. Assume that 
$T:=\sum_{j\in I} \cH_j$ is bounded in $L^2$ with operator norm $A_1$, and suppose that the inequalities (2.7.1) hold
 for $(n_1,n_2)=(0,0)$, $(0,1)$ or $(1,0)$.
Then $T$  maps $H^1_{prd}$ to $L^1$ 
with norm $\le C(A_0+A_1)$ (in particular $C$  is independent 
of $I$).
\endproclaim

\demo{Proof} This is a straightforward consequence of a theorem by 
R. Fefferman \cite{10}  which says that it suffices to check the 
operator on rectangle atoms. 
Suppose that $f$ is supported on a rectangle
parallel to the coordinate axes, with center $(c_1,c_2)$ 
 with sidelength $2^{-\ell_1}\times 2^{-\ell_2}$ and that
 $\|f\|_2\le 2^{(\ell_1+\ell_2)/2}$, moreover  $f$ satisfies the strong cancellation condition
$\int f(x_1,x_2) dx_1=0$ and
$\int f(x_1,x_2) dx_2=0$. Fefferman's theorem says that if 
$T$ is $L^2$ bounded and if the estimate
$$
\iint_{|x_1-c_1|\ge 2^{-\ell_1+n}} |Tf(x)| dx_1 dx_2+
\iint_{|x_2-c_2|\ge 2^{-\ell_2+n}} |Tf(x)| dx_1 dx_2
\lc 2^{-n\eps}
\tag 2.12
$$
holds for some $\eps>0$ then $T$ maps $H^1_{prd}(\Bbb R^2)$  boundedly to $L^1(\Bbb R^2)$.
 Since we assume $L^2$ boundedness 
it suffices  to prove (2.12).

 We estimate the corresponding 
integrals for $T$ replaced with $\cH_j$.
We use the size estimate in 
(2.7.1)  obtaining the bound $O(2^{\ell_1-n-j})$ for the 
$L^1$ norm in $\{x:|x_1-c_1|\ge 2^{-\ell_1+n}\}$ and we use the cancellation in $y_1$ together with the estimate (2.7.1) for the $y_1$-derivative to get the bound $2^{-\ell_1+j}$. 
Thus 
$$
\iint_{|x_1-c_1|\ge 2^{-\ell_1+n}} |\cH_jf(x)| dx_1 dx_2\lc 
\min\{ 2^{-\ell_1+j}, 2^{\ell_1-j-n}\}.
$$
We sum in $j$ and estimate the first term on the left of (2.12) by 
$C2^{-n/2}$.

Similarly (using now cancellation with respect to the $y_2$ variable) 
we obtain
$$
\iint_{|x_2-c_2|\ge 2^{-\ell_2+n}} |\cH_jf(x)| dx_1 dx_2\lc 
\min\{ 2^{-\ell_2+a_j}, 2^{\ell_2-a_j-n}\}.
$$
Clearly the right hand side is $O(2^{-n/2})$.
Let $j_0$ be the maximal $j$ with $2^{2a_j}\le 2^{-n+2\ell_2}$. Then  there is an absolute  constant $C_1$ so that
$2^{a_j}\lc 2^{a_{j_0}}2^{j-j_0}$ if $j\le j_0-C_1$. Thus 
$\sum_{j\le j_0-C_1} 2^{-\ell_2+a_j}\lc 2^{-\ell_2+a_{j_0}}\lc 2^{-n/2}$.
Similarly if $j_1$ denotes the minimal $j$ with 
$2^{2a_j}\ge 2^{-n+2\ell_2}$ then there is $C_2$ so that  for $j\ge j_1+C_2$ we we have $2^{-a_j}\lc 2^{ -a_{j_1}}2^{j_1-j}$ and  thus 
$\sum_{j\ge j_1+C_2} 2^{\ell_2-a_j-n}\lc 
2^{\ell_2-a_{j_1}-n}\lc 2^{-n/2}$.
We have only a bounded number of terms with
$j_0-C_1\le j\le j_1+C_2$; for those we use the bound $O(2^{-n/2})$. 
Combining the three estimates yields the bound $2^{-n/2}$ for the second term in (2.12).\qed

\enddemo

\head{\bf 3. $\boldkey L^{\boldkey p}$-boundedness of the Fourier integral contributions}\endhead
We now give an outline of the proof of Theorem B and consider first the maximal operator.
In view of 
 Lemma 2.1 and Corollary 2.3 it suffices to consider the maximal function 
$$\sup_j|\sum_{k>a_j} \cL^k\cR^k_j\cL^k f|.\tag 3.1 $$
where the sup is extended over a finite index set $J$. 
We use a familiar square-function technique and dominate
$$\align
\sup_j|\sum_{k>a_j} &\cL^k\cR^k_j\cL^k f|\\
&\le 
\Big(\sum_j
\big|\sum_{a_j<k\le b_j} \cL^k\cR^k_j\cL^k f\big |^2\Big)^{1/2}+
\Big(\sum_j
\big|\sum_{k>b_j} \cL^k\cR^k_j\cL^k f\big |^2\Big)^{1/2} 
\tag 3.2
\endalign
$$


We define an operator  $\fM$  acting on $F\in L^p(\ell^p)$
and an operator 
 $\widetilde \fM$ 
acting on  $G\in L^p(\ell^q(\ell^2))$ by
 by
$$\aligned
(\fM F)_j&=\cM_j F_j
\\
(\widetilde \fM G)_{j,k}&=\cM_j G_{j,k};
\endaligned 
\tag 3.3
$$ here the $\ell^2$ norm is taken with respect to the $k$ variable. 
We denote by 
$\|\fM\|_{p,q}$ the 
$L^p(\ell^q)\to L^p(\ell^q)$ operator norm of $\fM$ and by
$\|\widetilde \fM\|_{p,q,2}$ the 
$L^p(\ell^q(\ell^2))\to L^p(\ell^q(\ell^2))$ operator norm of $\widetilde \fM$.

We follow M. Christ \cite{7} (see also Nagel, Stein and Wainger \cite{17}
for a closely related earlier argument) and observe

\proclaim{Lemma 3.1}  For $1\le p\le 2$
$$\align
&\|\fM\|_{p,2} \lc (1+\|\cM\|_{L^p\to L^p})^{1-p/2}
\tag 3.4
\\
&\|\widetilde \fM\|_{p,2,2}
 \lc (1+\|\cM\|_{L^p\to L^p})^{1-p/2}
\tag 3.5
\endalign
$$\endproclaim

\demo{\bf Proof}
Since the operators $\cM_j$ are bounded on $L^p$, uniformly in $j$, it is clear  that
the operator norm of $\fM$ on $L^p(\ell^p)$ is $O(1)$; the same applies to the vector-valued setting by which we see that
the operator norm of $\widetilde \fM$ on $L^p(\ell^p(\ell^2))$ is finite.

Since $|\cM_j f_j|\lc \cM [\sup_\nu|f_\nu|]$ we see that the operator norm  
of $\fM $ on $L^p(\ell^\infty)$ 
and the operator norm of $\widetilde  \fM$ on $L^p(\ell^\infty(\ell^2))$
are  bounded by the $L^p$ norm of $\cM$. Interpolation gives the assertion.\qed
\enddemo

We now consider the first term on the right hand side of  (3.2).
First observe that there is a pointwise bound 
$$
|\cR^k_j g|\lc M_{str}(\cM_j (|g|))(x).
\tag 3.6
$$
To see this we use integration by parts in $\tau$ to estimate 
$$|\cR^k_j g|\lc \int \frac{2^k\widetilde \phi_j(x_1-y_1) \widetilde \chi(x,y)}
{(1+2^k|s|)^{N}}
|g(y_1, \gamma(x_1,x_2,x_1-y_1)+s)| dy_1 ds
$$
and change variables  $s=\ga(x_1, x_2+u,x_1-y_1)-\gamma(x_1,x_2,x_1-y_1)$ which is  legitimate 
since $\ga_{x_2}$ is close to $1$.

 By Littlewood-Paley theory for the
operators $\cL_k$, the pointwise bound (3.6)
 and the Fefferman-Stein theorem for the strong maximal function we get
$$\align
&\Big\|
\Big(\sum_j
\big|\sum_{a_j<k\le b_j} \cL^k\cR^k_j\cL^k f\big |^2\Big)^{1/2}\Big\|_p
\lc\Big\|
\Big(\sum\Sb j,k:\\ a_j<k\le b_j\endSb
\big| M_{str}(\cM_j (|\cL^k f|))\big |^2\Big)^{1/2}\Big\|_p
\\&\lc\Big\|
\Big(\sum\Sb j,k:\\ a_j<k\le b_j\endSb
\big| \cM_j (\cL^k f)\big |^2\Big)^{1/2}\Big\|_p
\lc \|\widetilde \fM\|_{p,2,2}\Big\|
\Big(\sum\Sb j,k:\\ a_j<k\le b_j\endSb
\big| \cL^k f |^2\Big)^{1/2}\Big\|_p
\lc \|\widetilde \fM\|_{p,2,2}
\|f\|_p
\tag 3.7
\endalign
$$
where for the last application of Littlewood-Paley theory we have used that for fixed $k$
  the cardinality
of the set $\{j: a_j<k\le b_j\}$ is bounded.

Similar but somewhat 
 more complicated arguments apply  to the second term in (3.2).
We need to introduce additional  dyadic decomposition in the variable dual to $x_1$ and 
 define operators $\cP_l$,  $\cQ_l$, $\Pi_m$ by
$$
\align
\widehat {\cP_l f}(\xi) &= \beta_0(2^{-l}\xi_1)\widehat f(\xi)
\\
\widehat {\cQ_l f}(\xi) &= \beta_l(\xi_1)\widehat f(\xi)
\endalign
$$
and  $$\Pi_m = \cP_{m-a_j+10}-\cP_{m-b_j-10}$$ and decompose  
 for fixed $k$ 
the identity operator as
$$I=\cP_{j+k-b_j-10}+\Pi_{j+k}+ \sum_{l>j+k-a_j+10} \cQ_l.
\tag 3.8$$

Then  we change  variables $k=b_j+n$, $l=j+k-a_j+m=j+(b_j-a_j)+m+n$
and see that
$$\align
\Big\|
\Big(\sum_j
\big|\sum_{k>b_j} \cL^k\cR^k_j\cL^k f\big |^2\Big)^{1/2} 
\Big\|
&\lc 
\sum_{n>0}\Big\|
\Big(\sum_j\big|\cL^{b_j+n}\cR^{b_j+n}_j 
\cP_{j+n-10} \cL^{b_j+n} f\big|^2\Big)^{1/2}\Big\|_p
\\
&+
\Big\|\Big(\sum_j \big|\sum_{k>b_j} \cL^k\cR^k_j\cL^k\Pi_{j+k}  f\big |^2\Big)^{1/2}\Big\|_p
\\&+
\sum_{n> 0}\sum_{m>0}
\Big\|
\Big(\sum_j\big|\cL^{b_j+n}\cR^{b_j+n}_j 
\cQ_{j+n+m+b_j-a_j} \cL^{b_j+n} f\big|^2\Big)^{1/2}\Big\|_p.
\tag 3.9
\endalign
$$

We need to show part (i) of the 
following proposition (part (ii) will be needed
for the singular Radon transform).

\proclaim{Proposition 3.2} Let $p_0>1$, let  $p_0\le p\le  2$ and define
$\theta\in [0,1]$ by $(1/p_0-1/p)=\theta(1/p_0-1/2)$. Then for $n>0$,
$m>0$

(i) 
$$
\align
&\Big\|
\Big(\sum_j\big|\cR^{b_j+n}_j 
\cP_{j+n-10} \cL^{b_j+n} f\big|^2\Big)^{1/2}\Big\|_p
\lc 2^{- \theta n/2} \|\fM\|_{p,2}^{1-\theta} \|f\|_p
\tag 3.10
\\
&\Big\|\Big(\sum_j \big|\sum_{k>b_j} \cL^k \cR^k_j\cL^k\Pi_{j+k}  f\big |^2\Big)^{1/2}\Big\|_p
\lc  \|\widetilde \fM\|_{p,2,2}^{1-\theta} \|f\|_p
\tag 3.11
\\
&\Big\|
\Big(\sum_j\big|\cR^{b_j+n}_j 
\cQ_{j+n+m+b_j-a_j} \cL^{b_j+n} f\big|^2\Big)^{1/2}\Big\|_p
\lc 2^{- \theta( n+m)/2} \|\fM\|_{p,2}^{1-\theta} \|f\|_p
\tag 3.12
\endalign
$$

(ii)
$$
\align
&\Big\|
\Big(\sum_j\big|
\cP_{j+n-10} \cR^{b_j+n}_j \cL^{b_j+n} f\big|^2\Big)^{1/2}\Big\|_p
\lc 2^{- \theta n/2} \|\fM\|_{p,2}^{1-\theta} \|f\|_p
\tag 3.13
\\
&\Big\|
\Big(\sum_j\big|
\cQ_{j+n+m+b_j-a_j}\cR^{b_j+n}_j  \cL^{b_j+n} f\big|^2\Big)^{1/2}\Big\|_p
\lc 2^{- \theta( n+m)/2} \|\fM\|_{p,2}^{1-\theta} \|f\|_p
\tag 3.14
\endalign
$$
\endproclaim

If we use 
the Fefferman-Stein theorem for
 vector-valued maximal functions we see that (3.9) and Proposition 3.2 imply the
bound
$$\align\Big\|
\Big(\sum_j
\big|\sum_{k>b_j} \cL^k\cR^k_j\cL^k f\big |^2\Big)^{1/2}
\Big\|_p &\lc (1+\|\fM\|_{p,2}+\|\widetilde \fM\|_{p,2,2})\|f\|_p 
\\
&\lc (1+\|\cM\|_{L^p\to L^p})^{1-p/2}\|f\|_p
\tag 3.15
\endalign$$
where for the last inequality we have used Lemma 3.1.
This  in conjunction  with (2.5), (2.8) and (3.7)   shows  that
$$\|\cM\|_{L^p\to L^p} \lc (1+\|\cM\|_{L^p\to L^p}^{1-p/2})
\tag 3.16$$
which  implies of course the $L^p$ boundedness of $\cM$ for $1<p\le 2$, with bound independent of $J$. Since $\cM$ is bounded on $L^\infty$ the $L^p$ boundedness for $2<p<\infty$ follows as 
well. By the monotone convergence 
theorem this shows the
$L^p$ boundedness of the maximal operator in Theorem B.

We turn to the proof of  Proposition 3.2.
The main technical Lemma used here concerns $L^2$   estimates 
for the kernels $\cR^k_j\cP_{j+k-b_j-10}$
and
$\cR^k_j\cQ_l$.

\proclaim{Lemma 3.3}
We have for $n>0$, $m>0$
  $$ \|\cR^k_j\cQ_l\|_{L^2\to L^2}
 \lc 2^{-n/2}, 
\quad \text{ if } k=b_j+n, l=j+k-a_j+m,
\tag 3.17$$
and 
$$ \|\cR^k_j\cP_{j+k-b_j-10}\|_{L^2\to L^2} \lc 2^{-(n+m)/2}  
\quad \text{ if } k=b_j+n. 
 \tag 3.18 $$

 The estimates (3.17-18) also hold with $\cR^k_j$ replaced by its adjoint $(\cR^k_j)^*$.

\endproclaim

The proof of Lemma 3.3 will be given  in the next section. The estimates involving the adjoint operator are  only needed for estimating  the 
singular Radon transform.
Taking the lemma  for granted we can now give the

\demo{\bf Proof of Proposition 3.2}
The scheme of the proof is the same as for the chain of inequalities in
(3.7).
For the main term (3.11) we use the Littlewood-Paley inequality
$$\Big\|\Big(\sum_{j,k} 
|\cL^k\Pi_{j+k} f|^2\Big)^{1/2}\Big\|_p\lc \|f\|_p,
\quad 1<p<\infty,
\tag 3.19
$$
and with (3.19) the proof of (3.11) follows by the argument given 
in (3.7). The inequality (3.19) in turn follows  by the usual argument 
involving the 
Marcinkiewicz multiplier theorem and Rademacher functions 
(see \cite{26}).  Here it is necessary to show the $L^p$  
boundedness of the operators
$\sum_{j,k}\pm \cL^k\Pi_{j+k}$ (for any choice of $\pm$) 
and the doubling assumption on $g$ is crucially used here.
We note that (3.19)  is essentially  a version of the 
angular Littlewood-Paley theorem used in \cite{17}, \cite{5}, 
\cite{7}, \cite{25}  and elsewhere.

For the terms  (3.10), (3.12) we use that  for fixed $n$ the 
$L^p(\ell^2)$ norm of $\{\cL^{b_j+n} f\}_{j\in \bbZ}$ is bounded by $C\|f\|_p$ and the 
argument in (3.7) shows that the left hand sides of both (3.10) and (3.12)
 are dominated by $C\|\fM\|_{p,2}\|f\|_p$ if $p>1$. 

For $p=2$ we have 
better bounds by  Lemma 3.3; indeed the left hand side of (3.10) for $p=2$
is dominated by $C2^{-\eps n}\|f\|_2$, the left hand side of (3.11) by $C\|f\|_2$  and for (3.12) we obtain the bound
$C 2^{-\eps(n+m)}\|f\|_2$. Interpolation yields (3.10), (3.11) (3.12). The proof of (3.13) and (3.14) is analogous if we take part (iii) of Lemma 3.3 into account. \qed
\enddemo

\subheading{${\boldkey L}^{\boldkey p}$-boundedness of the singular Radon transform}
In view of Proposition 2.4 and estimate (2.5) 
 we have to prove the boundedness of the 
 Fourier integral operator $\cF$   given by 
$$\cF f=\sum_j \sum_{k>a_j} \cL^k \cR^k_j \cL^k f. \tag 3.20$$

By a Littlewood-Paley estimate in the $x_2$ variables we see that
$$\align
&\Big\|
\sum_j \sum_{a_j<k\le b_j} \cL^k \cR^k_j \cL^k f\Big\|_p
\\&\lc  
\Big\|
\Big(\sum_j
\big|\sum_{a_j<k\le b_j} \cL^k\cR^k_j\cL^k f\big |^2\Big)^{1/2}\Big\|_p
\\&\lc \|\widetilde \fM\|_{p,2,2}
\|f\|_p
\tag 3.21\endalign
$$
where the second inequality had already been shown in (3.7).
Next we use the decomposition (3.8) 
and obtain

$$
\Big\|\sum_j \sum_{k> b_j} \cL^k \cR^k_j \cL^k f\Big\|_p
\lc  \sum_{n\ge 10} I_n  \,+\,II + \sum_{n\ge 10}\sum_{m>0} III_{m,n}
$$
where 
$$\align
I_n &= \Big\|
\sum_j\cL^{b_j+n}
\cP_{j+n-10} \cR^{b_j+n}_j  \cL^{b_j+n} f\Big\|_p
\\II&=
\Big\|\sum_j\sum_{k>b_j} \cL^k \Pi_{j+k} \cR^k_j\cL^k f\|_p
\\
III_{m,n}&= \Big\|
\sum_j \cL^{b_j+n}
\cQ_{j+n+m+b_j-a_j} \cR^{b_j+n}_j  \cL^{b_j+n} f\Big\|_p.
\endalign
$$
Now using for fixed $n$ the Littlewood-Paley decomposition 
$\{\cL^{b_j+n}\}_{j\in \bbZ}$ we see that $I_n$ is estimated by the left hand side of (3.13) and thus 
by $C 2^{-n\eps(p)} \|f\|_p$ with $\eps(p)>0$. Similarly
$III_{m,n}$ is dominated  by 
$C 2^{-(m+n)\eps(p)} \|f\|_p$, by (3.14).

Finally by Littlewood-Paley theory 
$$II\lc
\Big\|\Big(\sum_j \sum_{k>b_j}|\cR_j^k \cL^k f|^2\Big)^{1/2}\Big\|_p:=
\widetilde{II}.
$$
Now we decompose as in (3.8), but to the right hand side of $\cR^j_k$.
Thus
$$\widetilde{II}\lc  \sum_{n\ge 10} I_n'  \,+\,II' + \sum_{n\ge 10}\sum_{m>0} III_{m,n}'
$$
where 
$I_n'$ is the left hand side of (3.10), and $III_{m,n}'$ is the 
left hand side 
of (3.12). 
Moreover 
$$
II'=\Big\|\Big(\sum_j\sum_{k>b_j}
|\cR^k_j \cL^k \Pi_{j+k} f|^2\Big)^{1/2}\Big\|_p$$
which is dominated by $C\|\widetilde \fM\|_{p,2}\|f\|_p$; here we
use again (3.19) and the
Fefferman-Stein inequality.

Since we have already established the $L^p$ bounds for the maximal 
operator we know now by (3.4), (3.5) 
that $\|\fM\|_{p,2}$, $\|\fM_{p,2,2}\|_{p,2,2}$ are $O(1)$ and thus the combination
of previous estimates shows the $L^p$ boundedness of the Fourier integral operator $\cF$ in (3.20). As pointed out above this 
yields the $L^p$ boundedness for the singular Radon transform, 
for $1<p\le 2$. The estimates can be applied to the adjoint operator 
which yields the estimates in  the complementary range $2<p<\infty$.\qed.

\head{\bf 4. Proof of Lemma 3.3}\endhead
The kernel
of
$\cR^k_j\cQ_l$ is given by
$$
K(x,y)=\iiint e^{i(\tau(y_2-\ga(x,x_1-z_1))+\xi_1(z_1-y_1))} 
\chi(x,z_1,y_2) \phi_j(x_1-z_1)
\beta_1 (2^{-k}\tau)\beta_1(2^{-l}\xi_1)
 d\xi_1d\tau dz_1.
\tag 4.1$$
By integration by parts with respect to $z_1$ we obtain
$$
K(x,y)=\iiint e^{i(\tau(y_2-\ga(x,x_1-z_1))+\xi_1(z_1-y_1))} 
i^{-1} a(x,z_1,y_2,\tau,\xi_1)
 d\xi_1d\tau dz_1
$$ where $a=a_1+a_2$ with
$$\align
a_1(x,z_1,y_2,\tau,\xi_1)&= 
\frac{\ddot \gamma(x,x_1-z_1)}{(\tau\dot\ga(x,x_1-z_1)+\xi_1)^2}
 \chi(x,z_1,y_2) \phi_j(x_1-z_1)
\beta_1 (2^{-k}\tau)\beta_1 (2^{-l}\xi_1)
\\
a_2(x,z_1,y_2,\tau,\xi_1)&= 
\frac{\partial_{z_1}\big( \chi(x,z_1,y_2) \phi_j(x_1-z_1)\big)}{\tau\dot\ga(x,x_1-z_1)+\xi_1}
\beta_1 (2^{-k}\tau)\beta_1 (2^{-l}\xi_1).
\endalign
$$
The localization of the symbol in (4.1) implies that here $\xi_1|\gg |\tau \dot \gamma(x,x_1-z_1)|$.

The following fact will be crucial in the estimation of the $L^1$ norms in $x$ or $y$.
\proclaim{Sublemma} For large $j$ 
we have the estimates

(i)  $$\int_{|z_1-x_1|\le 2^{-j+1}} 2^{-j}|\ddot \gamma(x,x_1-z_1)| dz_1\lc 2^{-a_j}, \tag 4.2$$
(ii) $$\int_{2^{-j-1}\le |z_1-x_1|\le 2^{-j+1}} 2^j\frac{|\ddot \gamma(x,x_1-z_1)|}{(\dot \gamma(x,x_1-z_1))^2}
dz_1\lc 2^{b_j}. \tag 4.3$$

\endproclaim
\demo{Proof}
Let $I_j(x_1):=\{ z_1:0\le x_1-z_1\le 2^{-j+1}\}$.
By the quasimonotonicity assumption
we have for $t>0$ that $\ddot \gamma(x,t)=a(x,t)+O(\dot \gamma(x,t))$ where $a$ does not change sign and thus
$$
\align & 
\int\Sb I_j(x_1)\endSb
|\ddot \gamma(x,x_1-z_1)| dz_1\lc 
\Big| \int\Sb I_j(x_1)\endSb
\ddot \gamma(x,x_1-z_1) dz_1\Big|
+
 \int\Sb I_j(x_1)\endSb
|\dot \gamma(x,x_1-z_1)| dz_1
\\&\lc g(A2^{-j-1})\lc 2^{j-a_j}.
\endalign
$$
The same bound holds for the contribution over $x_1\le z_1$. This proves 
(4.2).

Now we turn to (4.3) and let
$J_j(x_1)
:=\{ z_1:2^{-j-1}\le x_1-z_1\le 2^{-j+1}\}$. 
Again by the quasimonotonicity of $\dot \gamma$ we get
$$
\int\limits\Sb J_j(x_1)\endSb
 \frac{|\ddot \gamma(x,x_1-z_1)|}{(\dot \gamma(x,x_1-z_1))^2}
dz_1
\le\Big|\int\Sb J_j(x_1)\endSb
 \frac{\ddot \gamma(x,x_1-z_1)}{(\dot \gamma(x,x_1-z_1))^2}
dz_1\Big|
+
\int\Sb J_j(x_1)\endSb|\dot \gamma(x,x_1-z_1)|^{-1}
dz_1.
$$
The first term on the right hand side  equals
$$
\Big|\frac 1{\dot \gamma(x, 2^{-j+1})}-
\frac 1{\dot \gamma(x, 2^{-j-1})}\Big|
\lc  [g(2^{-j-1}/A)]^{-1} \lc 2^{b_j-j}
$$
and the second term is estimated
by
$$\int\Sb J_j(x_1)\endSb
 g(2^{-j-1}/A) dz_1
\lc 2^{b_j-2j}\lc 2^{b_j-j}.$$
The contributions for
$2^{-j-1}\le z_1-x_1\le 2^{-j+1}$ are estimated in the same way and (4.3) is proved.
\qed\enddemo

\demo{Proof of Lemma 3.3, cont}
Integration by parts yields
$$
|K(x,y)|\lc \iiint \frac{(I-2^{2l} \partial_{\xi_1}^2)^N
(I-2^{2k} \partial_{\tau}^2)^N
a(x,z_1,y_2,\tau,\xi_1)}
{(1+2^{2k}|y_2-\ga(x,x_1-z_1)|^2)^N
(1+2^{2l}|y_1-z_1|^2)^N} dz_1d\xi_1 d\tau.
$$

Assume that $|\cC|\le 1$,
 $|\tau|, |\xi_1|\ge 1$ and that either
$|\cC\tau|\ge 2|\xi_1|$ or $|\xi_1|\ge 2 |\cC\tau|$.
(actually, for the  present proof of (3.17) we  need this for 
$|\xi_1|\ge 2|\Cal C\tau|$).

 It is easy to verify that  under these assumptions
 we have the product type symbol estimates
$$
\align
&\Big|\partial_{\xi_1}^{\alpha_1}\partial_\tau^{\alpha_2}
\big(
\beta_1 (2^{-k}\tau)\beta_1 (2^{-l}\xi_1) 
\tau(\cC\tau+\xi_1)^{-2}\big) \Big|\lc 2^{-l\alpha_1 }2^{-k \alpha_2 }  
2^k(|\cC\tau|+|\xi_1|)^{-2},
\\
&\Big|\partial_{\xi_1}^{\alpha_1}\partial_\tau^{\alpha_2}
\big(
\beta_1 (2^{-k}\tau)\beta_1 (2^{-l}\xi_1) 
(\cC\tau+\xi_1)^{-1}\big)\Big| \lc 2^{-l\alpha_1 }2^{-k \alpha_2 } (|\cC\tau|+|\xi_1|)^{-1}.
\endalign
$$
We apply this with $\cC= \dot \gamma(x,x_1-z_1)$ and see that 
$$
\align
&\big|\partial_{\xi_1}^{\alpha_1}\partial_\tau^{\alpha_2}
a_1(x,z_1,y_2,\tau,\xi_1)\big|
\lc  2^{-l\alpha_1 }2^{-k \alpha_2 }  \Big[
\frac{2^{k+j}}{|\cC\tau|+|\xi_1|)^{2}}
+\frac{2^{j}}{|\cC\tau|+|\xi_1|}\Big],
\\
&\big|\partial_{\xi_1}^{\alpha_1}\partial_\tau^{\alpha_2}
a_2(x,z_1,y_2,\tau,\xi_1)\big|
\lc  2^{-l\alpha_1 }2^{-k \alpha_2 }  
\frac{2^{2j}}{|\cC\tau|+|\xi_1|}.
\endalign
$$
Consequently
$$\multline
|K(x,y)|
\lc 
\iiint\limits \Sb|\xi_1|\approx 2^{l}
\\|\tau|\approx 2^k\\|x_1-z_1|\le 2^{-j+1}\endSb
\Big[\frac{2^{j+k} |\ddot \gamma(x,x_1-z_1)}
{(|\tau\dot \gamma(x,x_1-z_1)|+|\xi_1|)^2}
+ \frac{2^{2j}}{|\tau\dot \gamma(x,x_1-z_1)|+|\xi_1|}
\Big]\\
\times
(1+2^{2k}|y_2-\ga(x,x_1-z_1)|^2)^{-N}
(1+2^{2l}|y_1-z_1|^2)^{-N} dz_1d\xi_1 d\tau.
\endmultline
\tag 4.4
$$
We now examine the 
$L^1$ norm in $y$. We interchange the order of integration and first  integrate out in the $y$-variable. We take into account that now
$|\xi_1|\ge 2|\tau\dot\gamma(x,x_1-z_1)$ and obtain
$$\int|K(x,y)| dy
\lc 
\iiint\limits \Sb|\xi_1|\approx 2^{l}
\\|\tau|\approx 2^k\\|x_1-z_1|\le 2^{-j+1}\endSb
\big[2^{j+k}|\ddot \gamma(x,x_1-z_1)||\xi_1|^{-2}
+ 2^{2j}|\xi_1|^{-1}
\big] 2^{-k-l} dz_1d\xi_1 d\tau.
\tag 4.5
$$
By part (i) of  the Sublemma this is dominated by a constant times
$$
2^{-l-k}\iint\limits \Sb|\xi_1|\approx 2^{l}
\\|\tau|\approx 2^k\endSb [2^{2j+k-a_j-2l} +2^{j-l}] d\xi_1 d\tau \lc
\big[2^{2j+k-a_j-2l} +2^{j-l}\big]
\lc 2^{-m-n}.
\tag 4.6
$$
It is possible to show the same inequality for 
$\int|K(x,y)| dx$ but we can get away with the bound $O(1)$ for the latter integral and still get
(3.18).

We proceed similarly for the estimation of the kernel
 $\widetilde K$ of $\cR^k_j \cP_{j+k-b_j-10}$. Now however we have 
the restrictions $|\xi_1|\lc 2^{j+k-b_j-9}$ and 
$|\tau\dot\gamma(x,t)|\ge 2^{k-1}g(2^{-j-1}/A)\ge 2^{k-1+j-b_j}$
so that the latter expression is dominant.

The above analysis leads to

$$\multline
|\widetilde K(x,y)|
\lc 
\iiint\limits \Sb|\xi_1|\le 2^{j+k-b_j-9}
\\2^{k-1}\le |\tau| \le 2^{k+1}\\|x_1-z_1|\le 2^{-j+1}\endSb
\Big[\frac{2^{j+k} |\ddot \gamma(x,x_1-z_1)}
{|\tau\dot \gamma(x,x_1-z_1)|^2}
+ \frac{2^{2j}}{|\tau\dot \gamma(x,x_1-z_1)|}
\Big]\\
\times
(1+2^{2k}|y_2-\ga(x,x_1-z_1)|^2)^{-N}
(1+2^{2(j+k-b_j)}|y_1-z_1|^2)^{-N} dz_1d\xi_1 d\tau.
\endmultline
$$
and as above we get
$$
\int |\widetilde K(x,y)| dy
\lc  
\int\limits \Sb
|x_1-z_1|\le 2^{-j+1}\endSb
\Big[\frac{2^{j+k} |\ddot \gamma(x,x_1-z_1)|}
{|2^k\dot \gamma(x,x_1-z_1)|^2}
+ \frac{2^{2j}}{|\tau\dot \gamma(x,x_1-z_1)|}
\Big] dz_1d\xi_1 d\tau \lc 2^{b_j-k} \lc 2^{-n}
$$
where for the second to last inequality we use  (4.3).
Combining this with 
$\int |\widetilde K(x,y)| dx=O(1)$ we obtain (3.17).\qed
\enddemo

\enddocument